\documentclass[a4paper,fleqn,10pt]{article}
\usepackage{a4wide}
\usepackage[pdftex]{color}
\definecolor{urlcolor}{rgb}{0,0.5,0}
\definecolor{linkcolor}{rgb}{0.5,0,0}
\definecolor{citecolor}{rgb}{0,0,0.5}
\usepackage[pdfpagemode=UseOutlines,urlcolor=urlcolor,linkcolor=linkcolor,citecolor=citecolor,colorlinks=true]{hyperref}
\usepackage{amsmath}
\usepackage{amsfonts}
\usepackage{amssymb}
\usepackage[thmmarks,hyperref,amsmath]{ntheorem}
\newtheorem{Prop}{Proposition}[section]
\newtheorem{Thm}[Prop]{Theorem}
\newtheorem{Cor}[Prop]{Corollary}
\newtheorem{Lem}[Prop]{Lemma}
\theorembodyfont{\upshape}
\newtheorem{Rem}[Prop]{Remark}

\theoremstyle{nonumberplain}
\theoremheaderfont{\scshape}
\theorembodyfont{\upshape}
\theoremsymbol{\ensuremath{_\Box}}
\newtheorem{Proof}{Proof}
\usepackage[british]{babel}
\usepackage{enumerate}
\usepackage[pdftex]{graphicx}
\usepackage{natbib}
\bibpunct{(}{)}{;}{a}{}{,}

\DeclareMathOperator{\lspan}{span}

\DeclareMathOperator{\clos}{clos}
\DeclareMathOperator{\sinc}{sinc}
\newcommand{\N}{\mathbf{N}}

\newcommand{\Z}{\mathbf{Z}}
\newcommand{\R}{\mathbf{R}}
\newcommand{\C}{\mathbf{C}}
\DeclareMathOperator{\re}{\mathbf{Re}}

\renewcommand{\Re}{\re}

\newcommand{\rg}{\mathcal{R}}
\newcommand{\ke}{\mathcal{N}}

\newcommand{\cont}{\mathcal{C}}

\renewcommand{\L}{\mathrm{L}}

\newcommand{\st}{\ :\ }
\renewcommand{\H}{\mathcal{H}}
\newcommand{\HH}{\tilde\H}
\newcommand{\KK}{\tilde{K}}
\newcommand{\E}{\mathcal{E}}
\newcommand{\F}{\mathcal{F}}
\newcommand{\G}{\mathcal{G}}
\newcommand{\X}{\mathcal{X}}

\newcommand{\A}{\mathcal{A}}
\newcommand{\B}{\mathcal{B}}
\newcommand{\de}{\mathrm{d}}
\newcommand{\PW}{\mathrm{PW}}
\newcommand{\Ss}{\mathcal{S}}

\title{Representation by Integrating Reproducing Kernels}
\author{Thomas Hotz\footnote{Corresponding author: Thomas Hotz, Institute for Mathematical Stochastics, Georgia Augusta University of Goettingen, Goldschmidtstrasse 7, 37077 Goettingen, Germany; \texttt{hotz@math.uni-goettingen.de}.}~ and Fabian J. E. Telschow
\\[12pt]Institute for Mathematical Stochastics\\
Georgia Augusta University of Goettingen}

\begin{document}

\maketitle

\begin{abstract}
Based on direct integrals, a framework allowing to integrate a parametrised family of reproducing kernels with respect to some measure on the parameter space is developed. By pointwise integration, one obtains again a reproducing kernel whose corresponding Hilbert space is given as the image of the direct integral of the individual Hilbert spaces under the summation operator. This generalises the well-known results for finite sums of reproducing kernels; however, many more special cases are subsumed under this approach: so-called Mercer kernels obtained through series expansions; kernels generated by integral transforms; mixtures of positive definite functions; and in particular scale-mixtures of radial basis functions. This opens new vistas into known results, e.g. generalising the Kramer sampling theorem; it also offers interesting connections between measurements and integral transforms, e.g. allowing to apply the representer theorem in certain inverse problems, or bounding the pointwise error in the image domain when observing the pre-image under an integral transform.
\\[6pt]\noindent
\textbf{Keywords:} reproducing kernel; integral transform; radial basis function; scale-mixture; Kramer sampling; representer theorem.
\end{abstract}




\section{Overview}

Reproducing kernel Hilbert spaces (r.k.h.s.s) play an important r\^ole in many branches of mathematics, the corresponding reproducing kernels (r.k.s) also being called positive definite functions (p.d.f.s), radial basis functions (r.b.f.s) or autocorrelation functions (a.c.f.s) in special cases, see e.g. \citep{Wend:2005}. Often, r.k.s are constructed via integral transforms, cf. Section~\ref{sec:finite} and the references therein; by orthogonal expansions, so-called Mercer kernels, cf. Section~\ref{sec:mercer}; or by scale-mixtures of r.b.f.s, cf. Section~\ref{sec:rbf} and the references quoted there. In fact, all of these r.k.s share a common characteristic, as we will show below: they are obtained by integrating a parametrised family of r.k.s over the parameters; note that this also covers summation of kernels via integration with respect to a discrete measure.

Our aim is therefore to present an abstract framework, based on direct integrals, for the integration of r.k.s; this will be developed in Section~\ref{sec:framework} where Theorem~\ref{thm:int} clearly states conditions under which this is possible, allowing the r.k. to be calculated through pointwise integration, while also characterising its r.k.h.s. as the image of the direct integral of the individual Hilbert spaces under the summation operator. The aforementioned special cases of integral transforms, Mercer kernels, mixtures of p.d.f.s and r.b.f.s will be shown to be direct consequences of this theorem in Section~\ref{sec:special}, noting connections to some classical results of \cite{Bo:1933} and \cite{Schoe:1938a}. This framework also allows to view sampling equations, in particular Kramer sampling, from a slightly different and more general perspective in Section~\ref{sec:sampling}. Finally, we will point out some interesting relationships between the r.k.h.s. obtained by integration and the pre-images in the $\L_2$-space in Section~\ref{sec:data}, allowing the use of the representer theorem to solve an ``inverse problem'' over the $\L_2$-space in Proposition~\ref{prop:represent}, or bounding the pointwise error in the image domain by Proposition~\ref{prop:power}.

However, before we start to develop the abstract framework, we recall the basic definitions, and the simplest, classical case of the direct sum of two r.k.h.s. in Section~\ref{sec:sums}. This will form the starting point of which our abstract framework will be a generalisation.

\section{Direct sums of reproducing kernel Hilbert spaces}
\label{sec:sums}

Recall the notion of a \emph{reproducing kernel Hilbert space} (r.k.h.s.): let $\H$ be a Hilbert space of functions from some set $\X$ into $\C$ ---everything that follows applies equally if the scalar field is $\R$--- with scalar product $\langle \cdot, \cdot \rangle$; then the function $K$ defined on $\X \times \X$ is called a \emph{reproducing kernel} (r.k.) if for every $y \in \X$ we have $K(\cdot, y) \in \H$ as a function of the first argument, and $K$ possesses the \emph{reproducing property}, i.e. for all $f \in \H$ and $y \in \X$ we have
\begin{equation}
f(y) = \langle f, K(\cdot, y) \rangle \, .
\end{equation}
Then, $\H$ is called a r.k.h.s. over $\X$ with r.k. $K$.

In his seminal article, \cite{Aro:1950} proved that for two r.k.h.s. $\H_1$ and $\H_2$ over the same set $\X$ with scalar products $\langle \cdot, \cdot \rangle_1$ and $\langle \cdot, \cdot \rangle_2$, together with r.k. $K_1$ and $K_2$, respectively, their sum
\begin{equation}
\H = \H_1 + \H_2 = \{f = f_1 + f_2 \st f_1 \in \H_1, f_2 \in \H_2 \}
\end{equation}
is again a r.k.h.s., with r.k.
\begin{equation}
K = K_1 + K_2
\end{equation}
and norm $\Vert \cdot \Vert$ given by
\begin{equation}
\Vert f \Vert^2 = \inf_{f = f_1 + f_2 \st f_1 \in \H_1, f_2 \in \H_2} \Vert f_1 \Vert_1^2 + \Vert f_2 \Vert_2^2\,;
\end{equation}
see \citep[\S{}I.6]{Aro:1950}. He proved this by considering first the direct sum
\begin{equation}
\label{dirsum}
\HH = \H_1 \oplus \H_2
\end{equation}
with scalar product given by
\begin{equation}
\langle (f_1, f_2), (g_1, g_2) \rangle_\sim = \langle f_1, g_1 \rangle_1 + \langle f_2, g_2 \rangle_2\, .
\end{equation}
Then, $\H = \rg(S)$, the range of the summation operator 
\begin{equation}
\label{sumop}
S:\HH \rightarrow \H , \ (f, g) \mapsto f + g\, ,
\end{equation}
whose null space
\begin{equation}
\ke(S) = \{ (f, -f) \st f \in \H_1 \cap \H_2 \}
\end{equation}
is closed, whence $S$ is a bijection between $\ke(S)^\bot$ and $\rg(S) = \H$, and one can push forward the induced Hilbert space structure from $\HH$ to $\H$. The corresponding r.k. $K$ is then easily seen to be given by $K(\cdot, y) = S\KK(\cdot, y)$ where $\KK(\cdot, y) = (K_1(\cdot, y), K_2(\cdot, y))$. All one uses for this proof is the reproducing property of the respective kernels. It is clear that one can inductively obtain the r.k.h.s. corresponding to the sum of finitely many kernels.


\section{Abstract framework}
\label{sec:framework}

We now want to generalise the summation of kernels in the previous section to the integration of kernels. Towards this end, we first of all need an analogue of the direct sum in \eqref{dirsum}. For this, assume the index set $\Omega$ features a $\sigma$-algebra $\A$, and $\mu$ is a measure on that measurable space $(\Omega, \A)$. Moreover, for every $\omega \in \Omega$ let be given a r.k.h.s. $\H_\omega$ over $\X$ with scalar product $\langle \cdot, \cdot \rangle_\omega$, induced norm $\Vert \cdot \Vert_\omega$ and r.k. $K_\omega$.

These r.k.h.s.s need to be related in a measurable way: we assume that there is a partition of $\Omega$ into measurable sets $\Omega_n \in \A$, $n \in \N^{\infty} = \N \cup \{\infty\}$, and isometries $E_\omega: H_\omega \rightarrow \C^n$ for $\omega \in \Omega_n$ where we denote $\ell_2(\C)$ by $\C^\infty$ for uniformity in exposure. We then call a \emph{cross-section} $f = (f_\omega)_{\omega \in \Omega}$ of $f_\omega \in \H_\omega$, $\omega \in \Omega$ \emph{measurable} if for all $n \in \N^{\infty}$ the maps $\omega \mapsto E_\omega f_\omega : \Omega_n \rightarrow \C^n$ are measurable, whence $\omega \mapsto \Vert f_\omega \Vert_\omega : \Omega \rightarrow [0, \infty)$ is measurable, too.

A natural generalisation of the direct sum is then given by the \emph{direct integral}
\begin{align}
\HH &= \int^\oplus_{\omega \in \Omega} \H_\omega
\notag\\&= \Bigl\{ f = (f_\omega)_{\omega \in \Omega} \st f \text{ is a measurable cross-section and } \int_\Omega \Vert f_\omega \Vert_\omega^2 \de\mu(\omega) < \infty \Bigr\} \,,
\end{align}
which is again a Hilbert space with scalar product
\begin{equation}\label{eq:innerProd}
\langle f, g \rangle_\sim = \int_\Omega \langle f_\omega, g_\omega \rangle_\omega \, \de\mu(\omega)
\end{equation}
and norm $\Vert f \Vert_\sim = \langle f, f \rangle_\sim$ if we identify $f \in \HH$ and $g \in \HH$ in case $\int_\Omega \Vert f_\omega - g_\omega \Vert^2 \,\de\mu(\omega) = 0$. A worthwhile introdution into this topic can be found in \citep{Nie:2002}.

Note that the direct integral reduces to the direct sum if $\Omega$ is finite, $\A$ its power set, and $\mu$ is the counting measure on $\Omega$.

Next, we need a summation operator $S$ as in \eqref{sumop}, mapping $\HH$ into the space $\F(\X)$ of functions over $\X$ endowed with the topology of pointwise convergence. We define $S$ pointwise by setting
\begin{equation}
\label{eq:pointwiseIntegral}
(Sf)(x) = \int_\Omega f_\omega(x) \,\de\mu(\omega) = \int_\Omega \langle f_\omega, K_\omega(\cdot, x) \rangle_\omega \,\de\mu(\omega) = \sum_{n \in \N^\infty} \int_{\Omega_n} \langle E_\omega f_\omega, E_\omega K_\omega(\cdot, x) \rangle_{\C^n} \,\de\mu(\omega)
\end{equation}
for every cross-section $f \in \HH$ and every point $x \in \X$. The integrands in the last expression are clearly measurable if $\KK(\cdot, x) = (K_\omega(\cdot, x))_{\omega\in\Omega} \in \HH$, whence $(Sf)(x)$ is well-defined if
\begin{equation}
\int_\Omega \vert \langle f_\omega, K_\omega(\cdot, x) \rangle_\omega \vert \,\de\mu(\omega) < \infty.
\end{equation}
By applying Cauchy-Schwarz twice, we can estimate 
\begin{align}
\label{bound}
\int_\Omega \vert \langle f_\omega, K_\omega(\cdot, x) \rangle_\omega \vert \,\de\mu(\omega)
&\leq \int_\Omega \Vert f_\omega \Vert_\omega \, \Vert K_\omega(\cdot, x) \Vert_\omega \,\de\mu(\omega)
\notag\\&\leq \Biggl( \int_\Omega \Vert f_\omega \Vert_\omega^2 \,\de\mu(\omega) \Biggr)^\frac{1}{2} \Biggl( \int_\Omega \Vert K_\omega(\cdot, x) \Vert_\omega^2 \,\de\mu(\omega) \Biggr)^\frac{1}{2}
\notag\\&= \Vert f \Vert_\sim \, \Vert \KK(\cdot, x) \Vert_\sim
\notag\\&< \infty
\end{align}
by assumption. The cross-section $\KK(\cdot, x)$ is in $\HH$ for $x \in \X$ if it is measurable and its norm in $\HH$ is bounded; the latter is given by
\begin{equation}
\Vert \KK(\cdot, x) \Vert_\sim^2 = \int_\Omega \Vert K_\omega(\cdot, x) \Vert_\omega^2 \,\de\mu(\omega) = \int_\Omega K_\omega(x, x) \,\de\mu(\omega) \,.
\end{equation}

Furthermore, the operator $S:\HH \rightarrow \F(\X)$ is continuous; indeed, if a sequence $f^{(n)} \in \HH$ converges to some $f \in \H$ then
\begin{equation}
\Biggl\vert \int_{\Omega}f_{\omega}(x)-f_{\omega}^{(n)}(x)\de\mu(\omega) \Biggr\vert
\leq \int_{\Omega}\vert \langle f_{\omega}-f^{(n)}_{\omega},K_{\omega}(\cdot,x) \vert\de\mu(\omega)
\leq \Vert f-f^{(n)}\Vert_\sim \, \Vert \KK(\cdot, x) \Vert_\sim
\rightarrow 0
\end{equation}
for $n \rightarrow \infty$. Hence, the null space of $S$, 
\begin{equation}
\ke(S) = \Bigl\{ f \in \HH \st \int_\Omega f \,\de\mu(\omega) = 0 \Bigr\}\,,
\end{equation}
is closed and $S$ is a vector space isomorphism between $\ke(S)^\bot$ and $\H = \rg(S)$. We thus endow $\H$ with the Hilbert space structure turning $S$ into an isometry, i.e. we set for $f, g \in \HH$ the scalar product of $Sf, Sg \in \H$ to be
\begin{equation}
\label{pullback}
\langle Sf, Sg \rangle = \langle f, g \rangle_\sim \,.
\end{equation}
The obvious candidate for the r.k. on $\H$ is then given by
\begin{equation}
K(\cdot, x) = S\KK(\cdot, x)
\end{equation}
for every $x \in \X$.
Note that $\KK(\cdot, x) \in \ke(S)^\bot$ for every $x \in \X$ since for every $f \in \HH$
\begin{equation}
\label{eq:kernInNorth}
\langle f, \KK(\cdot, x) \rangle_\sim = \int_\Omega \langle f_\omega, K_\omega(\cdot, x) \rangle_\omega \, \de\mu(\omega) = \int_\Omega f_\omega(x) \,\de\mu(\omega) = (Sf)(x)\,,
\end{equation}
which is $0$ for $f \in \ke(S)$.
Furthermore, for any decomposition $f = f^{(1)} + f^{(2)} \in \HH$ with $f^{(1)} \in \ke(S)$, $f^{(2)} \in \ke(S)^\bot$ and any $x \in \X$
\begin{equation}
\langle Sf, K(\cdot, x) \rangle = \langle f^{(2)}, \KK(\cdot, x) \rangle_\sim = \langle f, \KK(\cdot, x) \rangle_\sim = (Sf)(x),
\end{equation}
which shows the reproducing property of $K$ for $\H$.

Finally, the norm $\Vert \cdot \Vert$ on $\H$ is given by
\begin{equation}
\Vert Sf \Vert^2 = \langle Sf, Sf \rangle = \langle f^{(2)}, f^{(2)} \rangle_\sim = \Vert f^{(2)} \Vert_\sim^2 =\inf_{g \in \HH \st Sg = Sf} \Vert g \Vert_\sim^2
\end{equation}
for $f = f^{(1)} + f^{(2)} \in \HH$ with $f^{(1)} \in \ke(S)$, $f^{(2)} \in \ke(S)^\bot$ since the map $f \mapsto f^{(2)}$ is the orthogonal projection onto $\ke(S)^\bot$.

Let us summarise what we have obtained:
\begin{Thm}
\label{thm:int}
Assume that for each index $\omega\in\Omega$ we are given a r.k.h.s. $\H_\omega$ over some set $\X$ with r.k. $K_\omega$ and norm $\Vert \cdot \Vert_\omega$; let $\HH$ be the direct integral of the $\H_\omega$ with respect to the measure space $(\Omega, \A, \mu)$. Furthermore assume that for every $x \in \X$ the cross-section of the r.k.s $\KK(\cdot, x) = (K_\omega(\cdot, x))_{\omega\in\Omega} \in \HH$, i.e. it is measurable and
\begin{equation}
\Vert \KK(\cdot, x) \Vert_\sim^2 = \int_{\Omega} K_{\omega}(x,x) \,\de\mu(\omega) < \infty.
\end{equation}
for every $x \in \X$.
Then,
\begin{equation}
\H = \Biggl\{ \int_\Omega f_\omega \,\de\mu(\omega) \st (f_\omega)_{\omega\in\Omega}\in \HH \Biggr\}
\end{equation}
is a r.k.h.s. over $\X$ with r.k. $K$ given by
\begin{equation}
K(x, y) = \int_\Omega K_\omega(x, y) \,\de\mu(\omega)
\end{equation}
for every $x,y \in \X$.
The norm $\Vert \cdot \Vert$ on $\H$ is given by
\begin{equation}
\Vert f \Vert^2 = \inf_{g \in \HH \st f = \int_\Omega g_\omega \,\de\mu(\omega)} \int_\Omega \Vert g_\omega \Vert^2 \,\de\mu(\omega)
\end{equation}
for any $f \in \H$. Also, for every $x \in \X$, $\KK(\cdot, x) \in \ke(S)^\bot$ with $\Vert \KK(\cdot, x) \Vert_\sim^2 = K(x,x)$.
\end{Thm}

\begin{Rem}
If $\mu$ is a finite measure, and if the r.k.s are uniformly bounded, i.e. there is a constant $c < \infty$ such that $K_\omega(x,x) < c$ for all $x \in \X$ and $\mu$-a.e. $\omega \in \Omega$, then $K(x,x) < c \mu(\Omega) < \infty$, so $K$ is uniformly bounded. If in addition, $\X$ is a topological space and for every $x \in \X$ and $\mu$-a.e. $\omega \in \Omega$ we have that $K_\omega(\cdot, x)$ is continuous, then $\H \subset \cont_b(\X)$, the latter denoting the space of uniformly bounded continuous functions over $\X$.
\end{Rem}

Let us note that the case of a direct sum of finitely many r.k.h.s.s considered in Section~\ref{sec:sums} is obtained as a special case of Theorem~\ref{thm:int} by endowing $\Omega = \{1, \dots, n\}$, $n \in \N$ with the counting measure.

A crucial -- but not easily verifiable -- ingredient to obtain a r.k.h.s. via integration is that $\KK(\cdot, x) = (K_\omega(\cdot, x))_{\omega\in\Omega}$ is a measurable cross-section in the sense we introduced at the beginning of this section, i.e. there exists a collection of isometries $\lbrace E_\omega: \H_\omega \rightarrow C^n \rbrace_{\omega\in\Omega}$ such that $\omega \mapsto E_{\omega} K_{\omega}(\cdot,x)$ is measurable. We will now state some conditions which guarantee this and are easy to check. The following Proposition and Lemma are taken from \citep[\S 2.8]{Nie:2002}.
\begin{Prop}\label{prop:constructIsom}
Let $\Omega$ be a measure space and $\lbrace H_{\omega}\rbrace_{\omega\in\Omega}$ Hilbert spaces. Suppose that $F$ is a countable set of cross-sections such that
\begin{enumerate}
\item for each $\omega\in\Omega$ the family $\lbrace f_\omega \rbrace_{f \in F}$ is dense in $H_{\omega}$; and
\item the map $\omega \mapsto \langle f_\omega, g_\omega \rangle_{\omega}$ is measurable for all $f,g \in F$.
\end{enumerate}
Then there exists a collection of isometries $\lbrace E_{\omega} \rbrace_{\omega\in\Omega}$ such that all $f\in F$ are measurable cross-sections. This collection of isometries is unique in the sense that if $\tilde E_{\omega}$ is another such collection then their sets of measurable cross-sections agree.
\end{Prop}
\begin{Lem}\label{lem:mbsections}
Let $F$ and $\lbrace E_{\omega} \rbrace_{\omega\in\Omega}$ be as in the Proposition \ref{prop:constructIsom}. Then a necessary and sufficient condition that a cross-section $g$ is measurable w.r.t. the collection $\lbrace E_{\omega}\rbrace_{\omega\in\Omega}$ is that $\omega \mapsto \langle g(\omega) , f(\omega) \rangle_\omega$ is measurable for all $f \in F$. 
\end{Lem}

Hence in the case of all $\H_{\omega}$ being r.k.h.s. over a separable topological space $\X$, we have the following corollaries:
\begin{Cor}
\label{cor:meascoh}
Let $\X$ be a separable topological space, $U \subset \X$ a countable, dense subset. Assume that $\{K_{\omega}(\cdot,y)\st y \in U\}$ is dense in $\H_{\omega}$ for each $\omega$ and that $\omega \mapsto K_{\omega}(x,y)$ is measurable for every $x\in\X$ and $y\in U$. Then there exists a collection of isometries $\lbrace E_{\omega} \rbrace_{\omega\in\Omega}$ such that $(K_\omega(\cdot,x))_{\omega\in\Omega}$ is a measurable cross-section for all $x\in\X$.
\end{Cor}

\begin{Cor}
\label{cor:contcoh}
Assume $\H_{\omega}$ for $\omega \in \Omega$ are r.k.h.s.s over a separable topological space $\X$ with r.k.s $K_{\omega}$. Moreover, assume that the maps
\begin{enumerate}
\item $x \mapsto K_\omega(x,x)$, as well as
\item $x \mapsto K_\omega(x,y)$ for all $y$ in a countable, dense subset $U$ of $\X$
\end{enumerate} 
are continuous for all $\omega \in \Omega$, and furthermore the maps $\omega \mapsto K_{\omega}(x,y)$ are measurable for all $x,y \in U$. Then there exists a collection of isometries $\lbrace E_{\omega} \rbrace_{\omega\in\Omega}$ such that $(K_\omega(\cdot,x))_{\omega\in\Omega}$ is a measurable cross-section for every $x\in\X$.
\end{Cor}
\begin{Proof}
Obviously, we want to use Proposition~\ref{prop:constructIsom}. Therefore we need to show that
\begin{equation}
\H_{\omega} = \clos_{\H_{\omega}}\left({\lspan\lbrace K_{\omega}(\cdot,y) \st y\in U \rbrace}\right).
\end{equation}
This is an immediate consequence of the continuity assumptions, the denseness of $U$, and $\H_{\omega} = \clos_{\H_{\omega}}\left({\lspan\lbrace K_{\omega}(\cdot,x) \st x\in \X \rbrace}\right)$; in fact, for every $\varepsilon > 0$ and $x \in \X$ there exists some $y \in U$ such that
\begin{equation}
\Vert K_{\omega}(\cdot,x) - K_{\omega}(\cdot,y) \Vert_{\omega}^2 = K_{\omega}(x,x) - 2\Re(K_{\omega}(x,y)) + K_{\omega}(y,y) <\epsilon.
\end{equation}
By Lemma~\ref{lem:mbsections} the map $\omega \mapsto K_{\omega}(\cdot, x)$ is measurable for $x\in \X \backslash U$ iff $\omega \mapsto K_{\omega}(x,y)$ is measurable for all $y \in U$. The latter is true since we can choose a sequence $x_n \in U$, $n\in\N$ with $\lim_{n\rightarrow\infty} x_n = x$, whence $\omega \mapsto K_{\omega}(x,y)$ is the limit of the measurable functions $\omega \mapsto K_{\omega}(x_n,y)$ and thus itself measurable.
\end{Proof}
\begin{Rem}
Note that Theorem~\ref{thm:int} can be viewed as a special case of the abstract framework in \citep{Sai:1983}: for some Hilbert space $\HH$ consider a map $k: \X \rightarrow \HH$. Then the set $\H = \{ x \mapsto \langle f, k(x) \rangle_\sim \st f \in \HH \}$ of functions over $\X$ forms a r.k.h.s. with kernel $K(x, y) = \langle k(y), k(x) \rangle_\sim$. For $\HH$ the direct integral above, $f$ any cross-section, and the special cross-sections $k(x) = (K_\omega(\cdot, x))_{\omega\in\Omega}$ this abstract construction in fact coincides with the one described above; in this setting, $S$ is given by $S(f)(x) = \langle f, k(x) \rangle_\sim$.
\end{Rem}
\begin{Rem}
Along the same lines as above, Theorem~\ref{thm:int} may be generalised to Hilbert space valued r.k.h.s.s (h.r.k.h.s.s). Recall the notion of a h.r.k.h.s., see e.g. \citep{CaDeViToi:2006}: in analogy to the scalar case, this is a Hilbert space $\H$ of functions over a set $\X$ with values in some separable Hilbert space $\E$ such that point evaluation is continuous. Denoting by $\B(\E)$ the set of bounded linear mappings from $\E$ to itself, continuity of point evaluation is, by the Riesz representer Theorem, indeed equivalent to the existence of a kernel $K:\X \times \X \rightarrow \B(\E)$, with the property that for all $x \in \X$ and $w \in \E$, the mapping $y \mapsto K(y,x)w$ is in $\H$, while fulfilling the reproducing equation
\begin{equation}
\langle f, K(\cdot,x)w \rangle_{\H} = \langle f(x), w \rangle_{\E}
\label{eq:RepPropVec}
\end{equation}
for every $f \in \H$. Note that this r.k. $K(x,x)$ is necessarily self-adjoint and satisfies
\begin{equation}
\Vert f(x) \Vert_{\E} \leq \sqrt{\Vert K(x,x) \Vert} \Vert f \Vert_{\H}.
\end{equation}
The latter shows that, given a collection of h.r.k.h.s. $\H_\omega$ with r.k.s $K_\omega$ over a measure space $(\Omega,\mu)$, any cross-section in the corresponding direct integral, which by definition is a square integrable function in the sense of the Bochner integral, is Bochner integrable if
\begin{equation}
\int_{\Omega} \Vert K_{\omega}(x,x)\Vert \de \mu(x) < \infty,
\end{equation}
cf. \eqref{bound}. Moreover, it is easily verified by the properties of the Bochner integral that the operator defined by
\begin{equation}
K(x,y):~w ~\mapsto~ \int_{\Omega} K_{\omega}(x, y)w\, \de \mu(\omega) \quad \text{for all } x,y \in \X
\label{}
\end{equation}
is in $\B(\E)$, and $K$ satisfies the properties of a r.k. for the image of the direct integral under the summation operator $S$ obtained by a pointwise application of the Bochner integral, cf. \eqref{eq:pointwiseIntegral}, the image being endowed with the pullback inner product under $S$, cf. \eqref{pullback}.

Note that this is still a special case of \citep[Proposition~20, p.~170]{Schwa:1964} if one views h.r.k.h.s.s as continuously embedded subspaces of the space $\F(\X)$ of all functions over $\X$ endowed with the topology of pointwise convergence.
\end{Rem}
\section{Special cases}
\label{sec:special}

\subsection{Integrating finite-dimensional reproducing kernel Hilbert spaces}
\label{sec:finite}

It is well-known that a huge class of r.k.h.s.s are describable as the range of some integral transform, take e.g. the Paley-Wiener spaces. We state some general conditions under which such kernels can be obtained via Theorem~\ref{thm:int}; cf. \citep{Sai:1997} for an extensive treatment of this topic.
\begin{Cor}\label{cor:trafokern}
In the setting of Theorem~\ref{thm:int},  assume there is a function $k:\X\times\Omega\rightarrow \C$ such that for every $x \in \X$, the function $\omega \mapsto k(x, \omega)$ is in $\L_2(\Omega,\mu)$. Then the r.k.h.s. $\H$ associated with the kernel
\begin{equation}
K(x,y)=\int_{\Omega}k(x,\omega)\,\overline{k(y,\omega)}\,\de\mu(\omega)
\end{equation}
is given by Theorem~\ref{thm:int}, i.e. $\H$ is the range of the operator $S:\L_2(\Omega)\rightarrow \F(\X)$,
\begin{equation}
(Sa)(y) = \int_{\Omega}a(\omega)k(y,\omega)\,\de\mu(\omega)\,.
\end{equation}
\end{Cor}
\begin{Proof}
Let $\H_\omega=\lspan\{k(\cdot,\omega)\}\subset\F(\X)$ endowed with the obvious inner product of the coefficients, i.e. $\Vert k(\cdot,\omega) \Vert_\omega = 1$, the isometries being given by $E_\omega:\H_\omega \rightarrow \C$, $k(\cdot,\omega) \mapsto 1$. Then, the r.k.s are given by
\begin{equation}
K_{\omega}(x,y) = k(x,\omega)\,\overline{k(y,\omega)}.
\end{equation}
Therefore, the cross-section $\KK(\cdot, y)$ is measurable iff $\omega \mapsto \overline{k(y,\omega)}$, which we assumed.
We then have
\begin{equation}
\HH=\{\omega\mapsto a(\omega)k(\cdot,\omega)\st a\in\L_2(\Omega,\mu)\}
\end{equation}
with
\begin{equation}
\Vert \omega\mapsto a(\omega)k(\cdot,\omega) \Vert_\sim = \Vert a \Vert_{\L_2(\Omega,\mu)}.
\end{equation}
At last we have for all $x \in \X$,
\begin{equation}
\int_{\Omega} K_{\omega}(x,x)\,\de\mu(\omega) = \int_{\Omega} \vert k_{\omega}(x) \vert^2\,\de\mu(\omega) = \Vert k_\cdot(x) \Vert_{\L_2(\Omega, \mu)}^2 < \infty
\end{equation}
by assumption. 
Hence, Theorem~\ref{thm:int} is applicable.
\end{Proof}
\begin{Rem}
Note that one can extend the above theorem to the case in which each $\H_\omega$ is the span of a finite number of such maps $k_i:\Omega\times\X\rightarrow\C$ for $i=1,\dots,d$ in the obvious way without changing the proof.
\end{Rem}

We will now revisit two classical examples.

\subsubsection{Paley-Wiener space}

Consider the continuous Fourier transform $F:\L_2(\R) \rightarrow \L_2(\R)$, $f \mapsto \int_\R f(x)\, e^{-2\pi i \cdot x}\, \de x$, which is an isometry, let $\Omega \subset \R$ be compact, further let $\Ss(\Omega) = \{f \in \L_2(\R) \st \forall \omega \notin \Omega\ f(\omega) = 0 \}$ be the closed sub-space of functions with support in $\Omega$, and denote the corresponding Paley-Wiener space of $\Omega$-band-limited functions by $\PW(\Omega) = F^{-1}(\Ss(\Omega))$. Note that $\Ss(\Omega)$ is isometrically isomorphic to $\L_2(\Omega)$ by restriction to $\Omega$.

Applying Corollary~\ref{cor:trafokern} to the function
\begin{equation}
k: \R \times \Omega \rightarrow \C,~ ~ (x,\omega) \mapsto e^{2 \pi i x \omega}\, ,
\label{}
\end{equation}
which obviously satisfies the assumption, and using the Lebesgue measure on $\Omega$ leads to the Hilbert space
\begin{equation}
\H=\left\{ f = \int_{\Omega} a(\omega) e^{2 \pi i x \omega} \,\de\omega \st a\in\L_2(\Omega)\right\}
\label{}
\end{equation}
with norm $\Vert f \Vert = \Vert a \Vert_{\L_2(\Omega)} = \Vert f \Vert_{\L_2(\R)}$. The latter is true by Parseval's identity, the former by the isometry of $\L_2(\Omega)$ and $\Ss(\Omega)$. Hence $\H$ is indeed the space $\PW(\Omega)$.
Moreover, we can conclude that
\begin{equation}
K(\cdot, y) = \int_\Omega K_\omega(\cdot, y) \,\de\mu(\omega) = \int_\Omega \exp(i\omega(\cdot-y)) \,\de\mu(\omega)
\end{equation}
is the r.k. for $\PW(\Omega)$, e.g. for $\Omega = [-\frac{1}{2}, \frac{1}{2}]$, we obtain the r.k. for $\PW([-\frac{1}{2}, \frac{1}{2}])$ as
\begin{align}
K(x, y) &= \int_{-\frac{1}{2}}^\frac{1}{2} \exp(2\pi i\omega(x-y)) \,\de\mu(\omega) = \frac{\sin(\pi (x-y))}{\pi (x-y)} = \sinc(\pi (x-y))\,.
\end{align}
Note that the latter is the well known kernel of the Paley-Wiener Space, see e.g. \citep[Theorem 10.12]{Wend:2005}.

\subsubsection{Global Sobolev kernels}
The global Sobolev space $W^m_{2}(\R^d)$ is the space of functions with weak derivative up to order $m$ being in $\L_2(\R^d)$. One endows this space with the inner
product
\begin{equation}
\langle f,g\rangle_{W^m_{2}(\R^d)}=\int_{\R^d}\left( 1+\Vert \omega\Vert^2 \right)^{m} F f(\omega)\overline{F g(\omega)}\de\omega \,,
\end{equation}
where $F:\L_2(\R) \rightarrow \L_2(\R)$ denotes the Fourier transform again. It is well known, see e.g. \citep[eq. (5.26)]{NaWa:1991} or again \citep[Theorem 10.12]{Wend:2005}, that this space possesses the r.k.
\begin{equation}
\label{eq:sobolevKernel}
K(x,y)=\int_{\R^{d}}\left( 1+\Vert \omega\Vert^2 \right)^{-m}e^{2\pi i(x-y)^t\omega}\de\omega
\end{equation}
if $d<2m$. This can be seen from Corollary~\ref{cor:trafokern} by considering the function
\begin{equation}
k:\R^d\times\R^d\rightarrow\C,~(x,\omega)\mapsto e^{2\pi ix^t\omega},
\end{equation}
which is uniformly bounded by 1, while $\omega\mapsto k(x,\omega)$ is measurable for each $x\in\R$ since it is continuous. Moreover, choose the measure $\mu$ on $\R^d$ s.t.
\begin{equation}
\de\mu(\omega)=\left( 1+\Vert\omega\Vert^2 \right)^{-m}\de\omega,
\end{equation}
denoting the Lebesgue measure by $\de\omega$. Surely, we have $\mu(\R^d)<\infty$ by $d<2m$ and hence $k(x,\cdot)\in\L_2(\R^d,\mu)$. Therefore we can apply Corollary~\ref{cor:trafokern} to obtain a r.k.h.s. $\H$ fulfilling
\begin{align}\notag
	\langle f,g\rangle_{\H}&=\int_{\R^d}a(\omega)\overline{b(\omega)}\de\mu(\omega)\\\notag
	&=\int_{\R^d}( 1+\Vert \omega\Vert^2)^mFF^{-1}\bigl( a(\omega)( 1+\Vert\omega\Vert^2 )^{-m} \bigr)\overline{FF^{-1}\bigl( b(\omega)( 1+\Vert\omega\Vert^2)^{-m} \bigr)}\de\omega\\
	&=\int_{\R^d}( 1+\Vert \omega\Vert^2)^mF(f)\overline{F(g)}\de\omega=\langle f,g\rangle_{W^{m}_{2}}
\end{align}
for $f=\int_{\R^d}a(\omega)e^{2\pi ix^t\omega}\de\mu(\omega)$, $g=\int_{\R^d}b(\omega)e^{2\pi ix^t\omega}\de\mu(\omega)$ and $a,b\in\L_2(\Omega,\mu)$.
Note that here, $\ke(S) = \{0\}$ by the injectivity of the Fourier transform, i.e. the representations in the direct integral are unique, so the scalar product is obtained by integrating over the individual scalar products.
Hence $\H=W^m_{2}(\R^d)$ with r.k. given in \eqref{eq:sobolevKernel}.

\subsubsection{Expansion kernels}
\label{sec:mercer}

Consider a collection $\varphi_n \in \F(\X)$, $n \in \N$, of linearly independent functions over $\X$; let $\lambda_n>0$ for $n\in\N$ and assume $\sum_{n\in\N} \lambda_n \vert \varphi_n(x) \vert^2 < \infty$ for all $x\in\X$. Clearly, Corollary~\ref{cor:trafokern} is applicable if we let $\Omega = \N$, $\mu(\{\omega\}) = \lambda_\omega$ and $k(\cdot, \omega) = \varphi_\omega$. The resulting Hilbert space is
\begin{equation}
\H = \biggl\{ \sum_{n \in \N} \gamma_n \varphi_n \st \sum_{n \in \N} \vert \gamma_n \vert^2 \lambda_n^{-1} < \infty \biggr\}.
\end{equation}
with kernel given by
\begin{equation}
K(x, y) = \sum_{n \in \N} \lambda_n \varphi_n(x) \overline{\varphi_n(y)}.
\end{equation}

By the linear independence of the ansatz functions $\varphi_n$, the kernel of $S : \L_2(\Omega) = \ell_2\bigl((\lambda_n)_{n\in\N} \bigr) = \{ (a_n)_{n \in \N} \st \sum_{n \in \N} \lambda_n \vert a_n \vert^2 < \infty \} \rightarrow \H$, $(a_n)_{n \in \N} \mapsto \sum_{n \in \N} \lambda_n a_n \varphi_n$ is trivial, $\ke(S) = \{ 0 \}$, whence $S$ is an isometry between $\H$ and the space of sequences square-summable w.r.t. the weights $\lambda_n$. Putting $\gamma_n = \lambda_n a_n $, we obtain an isometry to the space $\ell_2\bigl((\lambda_n^{-1})_{n\in\N} \bigr)$ such that $f \in \H$ iff $f = \sum_{n \in \N} \gamma_n \varphi_n$ for $(\gamma_n)_{n \in \N} \in \ell_2\bigl((\lambda_n^{-1})_{n\in\N} \bigr)$, the space of sequences square-summable w.r.t. the weights $\lambda_n^{-1}$. Hence, $(\sqrt{\lambda_n} \varphi_n)_{n \in \N}$ forms an orthonormal basis of $\H$.

\subsection{Integrating infinite-dimensional reproducing kernel Hilbert spaces}

We will now consider cases in which the individual r.k.h.s.s $\H_\omega$ are in general infinite-dimensional.

\subsubsection{Positive definite functions}

Assume that $\X$ is a group with neutral element $e$. We then call $\varphi: \X \rightarrow \C$ a \emph{positive definite function} (p.d.f.) if it gives rise to a r.k. $H$ via
\begin{equation}
H(x, y) = \varphi(xy^{-1}), \quad x,y \in \X;
\end{equation}
such a kernel is called \emph{translation invariant}.

The classical example is a finite-dimensional Euclidean vector space, $\X = \R^d$, $d \in \N$. The famous Theorem of \cite{Bo:1933} characterises the p.d.f.s in this case completely: $\psi$ is a p.d.f. if and only if there is a finite measure $\mu$ on $\R^d$ such that
\begin{equation}
\psi(x) = \int_{\R^d} \exp(i \omega^t x)\, \de\mu(\omega), \quad \forall x \in \X.
\end{equation}
Since $\varphi_\omega : \X \rightarrow \C$, $x \mapsto \exp(i \omega^t x)$ is a p.d.f. for every $\omega$, $K_\omega(x,y) = \exp(i\omega^t(x -y))$ having been used several times above, the easy ``if'' part of Bochner's Theorem is a direct consequence of Theorem~\ref{thm:int} with Corollary~\ref{cor:contcoh}. In fact, we have more generally:

\begin{Cor}
Let $\X$ be a separable topological group, $(\Omega,\mu)$ a measure space, and $\varphi:\X \times \Omega \rightarrow \C$ a map such that $\varphi(x,\cdot)$ is in $\L_2(\Omega,\mu)$ for every $x \in \X$ while $\varphi(\cdot, \omega)$ is a continuous p.d.f. for every $\omega \in \Omega$. Then
\begin{equation}
\psi(\cdot) = \int_\Omega \varphi(\cdot, \omega)\, \de\mu(\omega)
\end{equation}
is a p.d.f. again.
\end{Cor}

\subsubsection{Radial basis functions}
\label{sec:rbf}

We will now consider radial basis functions associated with r.k.h.s.s: we call a mapping $\varphi: [0, \infty) \rightarrow \C$ a \emph{radial basis function} (r.b.f.) on the metric space $(\X, \Delta)$ if it gives rise to a r.k. $H$ via
\begin{equation}
H(x, y) = \varphi(\Delta(x, y)), \quad x,y \in \X.
\end{equation}
For separable Hilbert spaces $\X$, \cite{Schoe:1938a} gave a characterisation of its r.b.f.s as \emph{scale mixtures}: $\psi$ is a r.b.f. if and only if there exists a finite measure on $[0, \infty)$ such that, $J_\alpha$ denoting the $\alpha$-th Bessel function of the first kind,
\begin{equation}
\psi(\delta) = \int_0^\infty \Gamma\bigl(\tfrac{d}{2}\bigr) \, \bigl(\tfrac{2}{\omega\delta}\bigr)^{\frac{d-2}{2}} \, J_{\frac{d-2}{2}}(\omega\delta)\, \de\mu(\omega), \quad \forall \delta \in [0, \infty),
\end{equation}
in case $\X$ is $d$-dimensional for some $d \in \N$, or
\begin{equation}
\psi(\delta) = \int_0^\infty \exp\bigl(-(\omega \delta)^2\bigr)\, \de\mu(\omega), \quad \forall \delta \in [0, \infty),
\end{equation}
in case $\X$ is infinite-dimensional. Again we can generalise the simpler ``if'' parts of Schoenberg's theorems using Theorem~\ref{thm:int} and Corollary~\ref{cor:contcoh}:

\begin{Cor}
Let $(\Omega,\mu)$ be a measure space, and $\varphi: [0, \infty) \times \Omega\rightarrow \C$ a map such that $\varphi(\delta,\cdot)$ is in $\L_2(\Omega,\mu)$ for every $\delta \in [0, \infty)$ while $\varphi(\cdot, \omega)$ is a continuous r.b.f. on the separable metric space $(\X, \Delta)$ for every $\omega \in \Omega$. Then for any finite positive measure $\mu$ on $\Omega$,
\begin{equation}
\psi(\cdot) = \int_\Omega \varphi(\cdot,\omega)\, \de\mu(\omega)
\end{equation}
is also a r.b.f. on $(\X, \Delta)$.
\end{Cor}

In fact, the last result is rather well-known for autocorrelation functions, see e.g. \citep[p.~355]{Yag:1987}, where it has been used to construct stationary and isotropic Gaussian random fields with certain desirable properties, e.g. through scale-mixtures of the Euclidean hat function by \cite{Gnei:1999}. More generally, scale-mixtures of compactly supported r.b.f. have been considered by \cite{Buh:1998} since they lead to numerically favourable band matrices.

\section{Sampling}
\label{sec:sampling}

The classical Kramer sampling theorem provides a method for obtaining orthogonal sampling theorems in the setting of integral transforms. We show that Kramer's sampling theorem \citep{Kra:1959} can be viewed as a statement about orthogonal bases of $\ke(S)^\bot$:
\begin{Prop}
\label{prop:sampling}
Assume we are in the setting of Theorem~\ref{thm:int}, its assumptions being fulfilled, with the r.k.h.s.s $\H_{\omega}$ all being of dimension $d \in \N^\infty$. Moreover, assume there is a sequence of $y_n \in \X, n\in\N$ such that $\{\omega \mapsto E_{\omega}K_{\omega}(\cdot,y_n)\}_{n\in\N}$ forms a complete orthogonal set of $\L_2(\Omega\rightarrow\C^d,\mu)$. Then $\{K(\cdot, y_n)\}_{n \in \N}$ forms a complete orthogonal set in $\H$; in particular one has
\begin{equation}
\label{eq:GenSampling}
g(x) = \sum_{n\in\N} g(y_n) \frac{K(x,y_{n})}{K(y_n, y_n)}\, ,
\end{equation}
for all $g\in\H$ and $x \in \X$.
\end{Prop}
\begin{Proof}
Let $f,g \in \HH\cap\ke(S)^\bot$ using the same notations as in Section \ref{sec:framework}.
The statement that $\{K(\cdot, y_n\}_{n \in \N}$ forms an orthogonal set is due to
\begin{equation}
\langle E_\cdot f_{\cdot} , E_\cdot g_{\cdot} \rangle_{\L^2(\Omega \rightarrow\C^d,\mu)} = \langle S(f) , S(g) \rangle_{\H},
\end{equation}
cf. \eqref{eq:innerProd}. Moreover, assume $0=\langle E_\cdot f_{\cdot} , E_\cdot K_{\cdot}(\cdot\cdot, y_n) \rangle_{\L^2(\Omega \rightarrow\C^d,\mu)} = \langle S(f) , K(\cdot,y_n) \rangle_{\H}$ for all $n \in \N$. By completeness we obtain $E_{\omega}f_{\omega} = 0$ $\mu$-a.e. and hence
\begin{equation}
S(f)(x) = \int_{\Omega}\langle E_\omega f_{\omega} , E_\omega K_{\omega}(\cdot,x) \rangle_{\omega} \de\mu(\omega) = 0
\label{}
\end{equation}
for all $x \in \X$ establishing the completeness of $\{K(\cdot, y_n\}_{n \in \N}$. By continuity of point evaluation in $\H$,  \eqref{eq:GenSampling} follows.
\end{Proof}
Note that we restricted ourselves to Hilbert spaces with identical dimensions merely for notational simplicity.

Using Corollary~\ref{cor:trafokern}, an immediate consequence of Proposition~\ref{prop:sampling} is the Kramer sampling theorem, see \citep{Jerri:1977}.
\begin{Cor}
Let $\mu$ be a measure on a space $\Omega$, and $k:\X\times\Omega\rightarrow \C$ a function such that for every $x \in \X$, the function $\omega \mapsto k(x, \omega)$ is measurable and in $\L_2(\Omega,\mu)$. If there exists a sequence of $y_n \in \X, n\in\N$ such that the set $\{k(y_n,\cdot)\}_{n\in\N}$ forms an orthogonal sequence in $\L^2(\Omega,\mu)$, then, for every $a\in\L_2(\Omega,\mu)$, we have the sampling equation
\begin{equation}
g(x)
 = \int_{\Omega} a(\omega) k(x, \omega) \de\mu(\omega)
 = \sum_{n\in\N} g(y_n) \frac{ \int_{\Omega} k(x,\omega)\overline{k(y_n, \omega)} \de\mu(\omega) }{ \int_{\Omega} \vert k(y_n, \omega) \vert^2 \de\mu(\omega) }\,.
\end{equation}
\end{Cor}

\begin{Rem}
Let us point out that in the specific case where $\X = \R^d$ and $y_n \in \Z^d$ for $n \in \N$, it is well known that the set $\{k(y_n,\cdot)\}_{n\in\N}$ forms an orthogonal sequence iff the bracket $[k(0, \cdot), k(0, \cdot)]$ is constant, the latter being defined through $[f, g](\omega) = \sum_{n \in \Z^d} (Ff)(\omega + n) \overline{(Fg)(\omega + n)}$ where $F$ again denotes the Fourier transform; see e.g. \citep{JePlo:2001} for details and extensions in this direction of shift invariant spaces.
\end{Rem}

\section{Measurements and representation}
\label{sec:data}

In the situation of Corollary~\ref{cor:trafokern}, we consider measurements either of the image $f = Sa$ of $a$ under $S$, or of $a$ itself. In the former case, we shall determine the corresponding pre-image $a \in \L_2(\Omega, \mu)$, in the latter case we shall be interested in determining the pointwise error in the image domain when interpolating the measurements of the pre-image $a$.

\subsection{Representation when observing the image}

First of all let us shortly recall the advantage of using a r.k.h.s. $\H$ for Tikhonov-like regularisation, which rests on the fact that every minimiser can be
expressed as a linear combination of the r.k. $K(\cdot,\cdot)$ evaluated at the sampling points $\{x_i\}_{i=1,\dots,N}$ due to the representer theorem. The following version of this theorem is due to \cite{Schoe:2000}:
\begin{Thm}
\label{thm:representer}
Let $\H$ be a r.k.h.s. $\H$ with r.k. $K:\X\times\X\rightarrow\C$. Furthermore, let $\{x_i\}_{i=1,\dots,N}\subseteq\X$ be a set of sampling points,
$\lambda:\R_{\geq0}\rightarrow\R$ a non-decreasing function and $L:\C^N\rightarrow\R \cup \{\infty\}$ an arbitrary loss function. Then the functional $J:\H\rightarrow \R \cup \{\infty\}$ given by
\begin{equation}
\label{lossfunctional}
    J(f) = L(f(x_1),\dots,f(x_{N})))+\lambda(\Vert{f}\Vert)
\end{equation}
possesses a minimiser $f$ of the form
\begin{equation}
\label{minimiser}
    f=\sum_{i=1}^N\alpha_i K(\cdot,x_i) \text{ with } \alpha_{1},\dots,\alpha_{N}\in\C.
\end{equation}
Furthermore, if $\lambda$ is strictly monotonically increasing, every minimiser is of form \eqref{minimiser}.
\end{Thm}
\begin{Proof}
Let $f$ be a minimiser of \eqref{lossfunctional}. Consider the orthogonal projection $f_{\parallel}$ of $f$ onto the finite dimensional subspace
$\lspan\{K(\cdot,x_i)\st 1\leq i\leq N\}$ and denote by $f_\bot$ the orthogonal part of $f$. Then by the reproducing property
\begin{equation}
    f(x_i)=(f,K(\cdot,x_i))=(f_{\parallel}+f_{\bot},K(\cdot,x_i))=(f_{\parallel},K(\cdot,x_i))=f_{\parallel}(x_i)
\end{equation}
and hence $L(f(x_1),\dots,f(x_N))=L(f_{\parallel}(x_1),\dots,f_{\parallel}(x_N))$.
On the other hand, we have
\begin{equation}
    \Vert f\Vert^2=\Vert f_\parallel\Vert^2+\Vert f_\bot\Vert^2\geq\Vert f_\parallel\Vert^2 \,,
\end{equation}
and therefore $f_{\parallel}$ is also a minimiser of \eqref{lossfunctional} since $\lambda$ is non-decreasing. In the case of a strictly increasing $\lambda$, we get
$\lambda(\Vert f\Vert)>\lambda(\Vert f_\parallel\Vert)$ if $\Vert f_\bot\Vert > 0$, and hence $f = f_\parallel$.
\end{Proof}

Now, using the notation of Section~\ref{sec:framework}, and assuming that $\H$ arises from integrating r.k.h.s.s as in Theorem~\ref{thm:int}, consider the functional $\tilde J:\HH \rightarrow \R \cup \{\infty\}$
\begin{align}\nonumber
	\tilde{J}(g) &= L\bigl((S(g))(x_1),\dots,(S(g))(x_N)\bigr)+\lambda(\Vert g \Vert_{\sim})\\
		&= L\bigl((S(g))(x_1),\dots,(S(g))(x_N)\bigr)+\lambda\left(\sqrt{\Vert g_{\parallel} \Vert_{\sim}^2+\Vert g_\bot \Vert_{\sim}^2}\right)
\label{DirIntFunctional}
\end{align}
with $g_\bot\in \ke(S)^\bot$ and $g_{\parallel}\in \ke(S)$ the respective orthogonal projections of $g \in \HH$. Then, if $\lambda$ is strictly monotonically increasing, the minimiser
$g$ of \eqref{DirIntFunctional} is an element of $\ke(S)^\bot$. Since $S(\ke(S)^\bot)\cong\H$ we have established that minimising \eqref{DirIntFunctional} is equivalent to minimising \eqref{lossfunctional}, the minimisers being related by $f = S(g) \in \H$.
The representer theorem above then yields that every minimiser of \eqref{DirIntFunctional} fulfils
\begin{equation}
	S(g)=\sum_{i=1}^{N} \alpha_i K(\cdot,x_i)=\int_{\Omega}\sum_{i=1}^{N}\alpha_i K_{\omega}(\cdot,x_i)\de\mu(\omega) \,,
\label{MinimiserHilbert}
\end{equation}
and thus, cf. \eqref{eq:kernInNorth},
\begin{equation}
	g=\sum_{i=1}^{N}\alpha_i \KK(\cdot,x_i)\in\ke(S)^\bot \,.
\label{MinimiserDirInt}
\end{equation}

This calculation becomes particularly interesting when $\H$ is obtained via an integral transform, i.e. in the situation of Corollary~\ref{cor:trafokern}; it then yields a method of estimating preimages of $f\in\H$ from measurements. In fact, we then always reconstruct the preimage from $\ke(S)^\bot$; let
\begin{equation}
	S:~\L_2(\Omega, \mu)\rightarrow S(\L_2(\Omega, \mu)),~ ~(S a)(x)=\int_{\Omega}a(\omega)k(\omega,x)\de\mu(\omega).
\end{equation}
Now, minimising $\tilde J$ in \eqref{DirIntFunctional} is equivalent to minimising $\bar J : \L_2(\Omega, \mu) \rightarrow \R \cup \{\infty\}$,
\begin{equation}
	\bar{J}(a) = L\bigl((S(a))(x_1),\dots,(S(a))(x_N)\bigr)+\lambda(\Vert a \Vert_{L_2(\Omega)})
\label{TrafoFunctional}
\end{equation}
for $a\in\L_2(\Omega, \mu)$, and by the derivation above its minimiser takes the form $a(\omega)=\sum_{i=1}^{N}\alpha_i \overline{k(\omega,x_i)}$. In summary, we obtain the following result:
\begin{Prop}
\label{prop:represent}
Let $S$ be an integral transform with kernel $k:\X \times \Omega \rightarrow \C$ satisfying $\omega \mapsto k(x,\omega) \in \L_{2}(\Omega,\mu)$ for every $x\in\X$, as in Corollary~\ref{cor:trafokern}, and assume $L,\lambda$ satisfy the assumptions of Theorem~\ref{thm:representer}. Then the functional $\bar J: \L_2(\Omega, \mu) \rightarrow \R \cup \{\infty\}$,
\begin{equation}
\bar{J}(a) = L\bigl((S(a))(x_1),\dots,(S(a))(x_N)\bigr)+\lambda(\Vert a \Vert_{L_2(\Omega)}),
\end{equation}
possesses a minimiser $a^*$ admitting the representation
\begin{equation}
\label{prop:tikhonovSpec}
a^*(\omega) = \sum_{i=1}^N \alpha_{i} k(\omega, x_i).
\end{equation}
Indeed, for any minimiser $a^*$ of $\bar J$, we have that $S(a^*)$ minimises $J$ in Theorem~\ref{thm:representer}, while for any minimiser $f \in \H$ of $J$ the unique pre-image $a^* \in \ke(S)^\bot$ with $S(a^*) = f$ minimises $\bar J$.

Furthermore, if $\lambda$ is strictly increasing, then any minimiser $a^*$ of $\bar{J}$ is of this form. In fact, then $a^* \in \ke(S)^\bot$.
\end{Prop}
\begin{Rem}
A typical situation where the minimisation of $\bar J$ \eqref{TrafoFunctional} occurs is in \emph{inverse problems}: one only observes the image $S(a)$ of the function $a \in \L_2(\Omega, \mu)$ of interest, usually with some noise; here, $S$ is an integral operator. One then wants to find a function $a^*$ which is close to the data as measured by $L$ but not too large in norm, as the perturbed data no longer lie in the range of $S$, i.e. in the r.k.h.s $\H$, whence the regularisation via $\lambda$. There is then a well-developed theory showing under which conditions the minimiser $a^*$ will be close to the true function $a$, see e.g. \citep{EnHaNeu:1996}.

As an example consider for some regularisation parameter $\gamma > 0$ and data $y = (y_i)_{i=1}^N \in \C^N$ the quadratic loss function
\begin{equation}
\bar J(a) = \sum_{i=1}^N \vert y_i - (S(a))(x_i) \vert^2 + \gamma \Vert a \Vert_{\L_2(\Omega, \mu)}^2 \,,
\end{equation}
to be minimised over $a \in \L_2(\Omega, \mu)$. 
Then, with the matrix $H = \bigl( K(x_i, x_j) \bigr)_{i,j=1}^N \in \C^{N \times N}$, the loss in dependence of $\alpha$ is given by
\begin{equation}
\Vert H \alpha - y \Vert^2 + \gamma\, \biggl\Vert \sum_{i=1}^N \alpha_{i} K(\cdot,x_i) \biggr\Vert^2_{\H}
= \alpha^*H^*H\alpha - \alpha^*H^*y - y^*H\alpha + y^*y + \gamma\alpha^* H\alpha.
\label{Loss}
\end{equation}
A short calculation shows that the minimising $\alpha$ solves the equation
\begin{equation}
	(H^*H + \gamma H)\alpha = H^*y \,,
\end{equation}
and thus the minimiser $a^*$ of the functional can be computed explicitly from \eqref{prop:tikhonovSpec}.
\end{Rem}

\subsection{Interpolating the pre-image}

We now change our viewpoint, assuming that we observe the pre-image $a$. For this to make sense, in addition to the assumptions of Corollary~\ref{cor:trafokern}, let $\G$ be a r.k.h.s. over $\Omega$ with kernel $G$ such that the diagonal $\omega \mapsto d(w) = G(w,w) = \Vert G(\cdot, \omega) \Vert_\G \in \L_2(\Omega, \mu)$; thence $G(\cdot, \omega) \in \L_2(\Omega, \mu)$ for every $\omega$, so $\G \subset \L_2(\Omega, \mu)$ with a continuous embedding whose norm is bounded by the $\L_2(\Omega, \mu)$-norm of the diagonal:
\begin{equation}
\Vert a \Vert_{\L_2(\Omega, \mu)} = \Vert \omega \mapsto \langle a(\omega), G(\cdot, \omega) \rangle_\G \Vert_{\L_2(\Omega, \mu)}
\leq \Vert a \Vert_\G \Vert d \Vert_{\L_2(\Omega, \mu)}
\end{equation}
for every $a \in \G$. Now, let $W \subset \Omega$ be a set on which we observe $a \in \G$, and denote by $P_W$ the corresponding power function, given by
\begin{equation}
P_W(\omega)^2 = G(\omega, \omega) - G_W(\omega, \omega),
\end{equation}
where $G_W$ is the kernel of the sub-space $\G_W$ generated by $\{G(\cdot, \omega) \st \omega \in W\}$; observe that $P_W \in \L_2(\Omega, \mu)$, too.

We are interested in estimating the pointwise error at $x \in \X$ made by approximating $a$ by the interpolant $a_W \in \G_W$ with $a_W(\omega) = a(\omega)$ for all $\omega \in W$, i.e. for $f = S(a)$ and $f_W = S(a_W)$ we estimate
\begin{align}
\vert f(x) - f_W(x) \vert &= \vert \langle S(a) - S(a_W), K(\cdot, x) \rangle_\H \vert
\notag\\ &\leq \vert \langle a - a_W, k(x, \cdot) \rangle_{\L_2(\Omega, \mu)} \vert
\notag\\ &\leq \Vert a - a_W \Vert_{\L_2(\Omega, \mu)} \ \Vert k(x, \cdot) \Vert_{\L_2(\Omega, \mu)}
\notag\\ &\leq \Vert a \Vert_\G \ \Vert P_W \Vert_{\L_2(\Omega, \mu)} \ \Vert k(x, \cdot) \Vert_{\L_2(\Omega, \mu)} \,;
\label{l2bound}
\end{align}
observe that the first inequality is in fact an equality if $\ke(S) = \{0\}$. Recall that $a_W$ can also be characterised as the function in $\G$ with minimal norm interpolating $a(\omega)$ at all $\omega \in W$.
\begin{Prop}
\label{prop:power}
Let $(\Omega,\mu)$ be a measure space, $W\subset \Omega$ a set of sample points, $\X$ a set and $\G\subset\L_{2}(\Omega,\mu)$ a continuously embedded r.k.h.s. 
Moreover, assume that we are given an integral transform $S$ with kernels $k:\X \times \Omega \rightarrow \C$ such that $S(\L_{2}(\Omega,\mu))$ is a r.k.h.s. as in Corollary~\ref{cor:trafokern}. Then for $a \in \G$ and all $x\in\X$, the pointwise difference between the image of the minimum-norm interpolator $a_W$, with corresponding power function $P_W$, and the image of $a$ at $x$ can be bounded by
\begin{equation}
\vert S(a)(x) - S(a_{W})(x) \vert \leq \Vert a \Vert_{\G}\ \Vert P_W \Vert_{\L_{2}(\Omega,\mu)}\ \Vert k(x, \cdot) \Vert_{\L_{2}(\Omega,\mu)}.
\end{equation}
\end{Prop}
Note that, in order to put this to practical use, one will have to be able to compute $S(G(\cdot, \omega))$ to obtain the image $f_W$ of the interpolant $a_W$ under $S$; indeed, if $W$ is finite, $a_W$ is given by $\sum_{\omega \in W} \alpha_\omega G(\cdot, \omega)$ for some $\alpha_\omega \in \C$.

\paragraph{Acknowledgements.}
We are indebted to Prof. Dr. Mikhail Gordin of the St. Petersburg Department of the V. A. Steklov Institute of Mathematics, Russian Academy of Sciences, for equipping us with the analytical tools required to address and connect the issues arising here, as well as for numerous helpful discussions. Furthermore, we thank Prof. em. Dr. Robert Schaback, University of G\"ottingen, for encouraging us to pursue this research, and for helping us to connect it to the many areas where this applies and relates to. Finally, we thank the German Research Foundation (DFG) for support via DFG SFB 803; the questions raised in that Collaborative Research Centre in fact initiated this research.

\bibliographystyle{elsart-harv}
\bibliography{kern.bib}

\end{document}